\documentclass[a4paper,11pt]{amsart}

\usepackage{graphicx}
\usepackage{mathptmx}
\usepackage{amsmath}
\usepackage{amssymb}
\usepackage{enumitem}
\usepackage{xcolor} 
\usepackage{pgfplots}

\newmuskip\pFqmuskip

\newcommand*\pFq[6][8]{%
  \begingroup % only local assignments
  \pFqmuskip=#1mu\relax
  \mathcode`=\string"8000
  \begingroup\lccode`\~=`\,
  \lowercase{\endgroup\let~}\pFqcomma
  F^{#2}_{#3}{\left(\genfrac..{0pt}{}{#4}{#5}\bigg|#6\right)}%
  \endgroup
}
\newcommand{\pFqcomma}{\mskip\pFqmuskip}

\newcommand{\ap}{a^{\dagger}}
\newcommand{\ba}{aa^{\dagger}}
\newcommand{\fa}{a^{\dagger}a}

\newtheorem{theorem}{Theorem}[section]

\begin{document}

\title[]{Generalized Bell polynomial operators arising from generalized normal ordering}

\author{Taekyun  Kim}
\address{Department of Mathematics, Kwangwoon University, Seoul 139-701, Republic of Korea}
\email{tkkim@kw.ac.kr}
\author{Dae San  Kim}
\address{Department of Mathematics, Sogang University, Seoul 121-742, Republic of Korea}
\email{dskim@sogang.ac.kr}


\begin{abstract}
This paper explores the deformed combinatorial structures arising from the generalized Heisenberg algebra $\text{GHA}_f$, which is characterized by an analytic function of the Hamiltonian $f(H)$ and governs systems with non-linear spectra. Moving beyond the classical Heisenberg-Weyl framework, we investigate the normal ordering of the generalized number operator $N_f^n = (\fa)^n$, which naturally introduces the generalized Stirling operators of the second kind ${n \brace k}_f(H)$. Using the vacuum eigenvalue $\alpha_0$ of the Hamiltonian, we define quantum operator factorials and a generalized quantum exponential function $e_f^x$. We explicitly construct the generalized coherent states $|z\rangle$ and derive several operator identities. Notably, we prove that the powers of quantum operators expand into generalized falling factorials, and that the coherent state expectation values $\langle z| N_f^n |z\rangle$ are explicitly given by the generalized Bell polynomial operators $\phi_{n,f}(|z|^2)(H)$.
\end{abstract}

\maketitle

\markboth{\centerline{\scriptsize Generalized Bell polynomial operators arising from generalized normal ordering}}
{\centerline{\scriptsize Taekyun Kim and Dae San Kim}}

\section{Introduction} 
The study of algebraic structures in quantum mechanics has long been a cornerstone for understanding physical systems, with the classical Heisenberg-Weyl algebra serving as the foundational paradigm. Generated by the creation operator $\ap$ and annihilation operator $a$ satisfying the canonical commutation relation $[a, \ap] = 1$, this algebra perfectly describes the quantum harmonic oscillator. Within this framework, combinatorial structures naturally emerge. A prime example is the normal ordering of powers of the number operator $N = \fa$, where the coefficients are precisely the Stirling numbers of the second kind. These numbers, alongside their closely related Bell polynomials, bridge the gap between quantum operator algebra and classical combinatorics, providing deep insights into the statistical properties of quantum states, such as coherent states. In recent decades, the search for novel physical phenomena and more flexible mathematical frameworks has led to various generalizations of this standard structure. Among these, the generalized Heisenberg algebra (denoted here as $\text{GHA}_f$) stands out. Generated by the set $\{a, \ap, H\}$ and governed by an analytic function of the Hamiltonian $f(H)$, the $\text{GHA}_f$ allows for the description of systems with non-linear spectra and deformed symmetries. When transitioning from the classical algebra to its generalized counterpart, a natural question arises: how do the underlying combinatorial structures deform? Specifically, how do the traditional Stirling numbers and Bell polynomials evolve when the standard number operator is replaced by the generalized number operator $N_f = \fa$? The aim of this paper is to systematically explore this connection, developing the quantum operator calculus that arises from the normal ordering within $\text{GHA}_f$ and applying it to generalized coherent states. \\
{\bf{Structure of the Paper}} \\
To establish a self-contained framework, the remainder of this paper is structured as follows:\\
$\indent \bullet$ In Section 2, we establish the mathematical preliminaries. We first review the classical definitions of the Stirling numbers of the second kind and Bell polynomials. We then remind the reader of the standard Heisenberg-Weyl algebra, detailing the actions of $a$, $\ap$, and $N$ on the number states $|n\rangle$, the action of $a$ on coherent states $|z\rangle$, and the classic normal ordering expansion of $N^n$. Following this, we transition to the generalized Heisenberg algebra $\text{GHA}_f$, defined by the commutation relations:$$aH = f(H)a, \quad H \ap = \ap f(H), \quad \text{and} \quad [a, \ap] = f(H) - H.$$ Within this generalized setting, we introduce the quantum operator functions $[k]_f(H) = f^k(H) - H$ for any integer $k$, and derive several algebraic relations involving these functions and powers of the creation and annihilation operators. We conclude the section by analyzing the actions of these operators on the generalized number states and recalling the normal ordering of $N_f^n = (\fa)^n$, which naturally gives rise to the generalized Stirling operators of the second kind, denoted as ${n \brace k}_f(H)$.\\
$\indent \bullet$ In Section 3, we present the main results of this work, focusing on the construction of generalized coherent states and their combinatorial properties. Let $\alpha_0$ denote the eigenvalue of the Hamiltonian $H$ corresponding to the vacuum state. Utilizing this, we introduce the quantum operator factorials $[n]_f(\alpha_0)!$ and define the generalized quantum operator exponential function:$$e_f^x = \sum_{n=0}^{\infty} \frac{x^n}{[n]_f(\alpha_0)!}.$$ Our primary findings are organized as follows:\\
$\indent \bullet$ Theorem 3.1 provides an explicit expression of the generalized coherent state $|z\rangle$ in terms of the number states.\\
$\indent \bullet$ Theorem 3.2 establishes an operator identity connecting the powers of the quantum operator to generalized falling factorial operators $\big([m]_f(\alpha_0)\big)_l$ via the generalized Stirling operators:$$\big([m]_f(\alpha_0)\big)^n = \sum_{l=0}^{n} {n \brace l}_f(\alpha_{m}) \big([m]_f(\alpha_0)\big)_l$$\\
$\indent \bullet$ Theorems 3.3 and 3.4 determine the expectation value of the $n$-th power of the generalized number operator in a coherent state, showing that $\langle z| N_f^n |z\rangle = \phi_{n,f}(|z|^2)(H)$, where $\phi_{n,f}(x)(H) = \sum_{l=0}^{n} {n \brace l}_f(H) x^l$ represents the generalized Bell polynomial operators arising from normal ordering.\\
$\indent \bullet$ Theorem 3.5 derives a closed-form summation formula for these generalized Bell polynomial operators, proving that:$$\phi_{k,f}\big(|z|^2\big)(\alpha_{0}+|z|^{2}) = e_f^{-|z|^2} \sum_{n=0}^{\infty} \frac{|z|^{2n}}{[n]_f(\alpha_0)!} \big([n]_f(\alpha_0)\big)^k.$$\\
$\indent \bullet$ In Section 4, we conclude the paper by summarizing the results in this paper.

\section{Preliminaries}

It is well known that the Stirling numbers of the second kind are defined by 
\begin{equation}
x^{n}=\sum_{k=0}^{n}{n \brace k}(x)_{k},\quad (n\ge 0),\quad (\mathrm{see}\ [1-8]), \label{1}
\end{equation}
where 
\begin{equation*}
(x)_{0}=1,\quad (x)_{n}=x(x-1)(x-2)\cdots(x-n+1),\ (n\ge 1). 
\end{equation*}
The Bell polynomials are given by 
\begin{equation}
\phi_{n}(x)=\sum_{k=0}^{n}{n \brace k}x^{k},\quad (n\ge 0),\quad (\mathrm{see}\ [3,4,6-8]). \label{2}
\end{equation}
When $x=1$, $\phi_{n}=\phi_{n}(1)$ are called the Bell numbers. Boson normal ordering is the process of rearranging a Stirling of creation $\ap$ and annihilation $a$ operators so that all creation operators are to the left of all annihilation operators. \par 
For a single boson mode, the fundamental equation governing the ordering process is the Heisenberg-Weyl algebra commutation relation: $\ba-\fa=1$. The operators $\ap$ and $a$ act on the number states in Fock space as
\begin{equation}
a|n\rangle=\sqrt{n}|n-1\rangle,\quad \ap|n\rangle=\sqrt{n+1}|n+1\rangle,\quad (\mathrm{see}\ [1,7]). \label{3}
\end{equation}
The number operator $N$ is defined by 
\begin{equation}
N|n\rangle=n|n\rangle,\quad (\mathrm{see}\ [9-13]).\label{4}
\end{equation}
By \eqref{3} and \eqref{4}, we get $N=\ap a$. The states are assumed to be orthonormal so that $\langle m|n\rangle=\delta_{m,n}$, where $\delta_{m,n}$ is the Kronecker delta. \par 
The second set of states of interest in Fock space is the coherent states $|z\rangle,\ z\in\mathbb{C}$. They are defined as the eigenstates of the annihilation operator 
\begin{equation}
a|z\rangle=z|z\rangle,\quad\mathrm{equivalently}\quad \langle z|\ap=\langle z|\overline{z}, \quad\textrm{and $\langle z|z\rangle=1$, $(z\in\mathbb{C})$.}\label{5}
\end{equation}
Katriel showed that the Stirling numbers of the second kind appear as the coefficients in the normal ordering of the $n$-th power of the number operator, namely 
\begin{equation*}
N^{n}=(\fa)^{n}=\sum_{k=0}^{n}{n \brace k}(\ap)^{k}a^{k},\quad (n\ge 0)\quad (\mathrm{see}\ [3,4]).
\end{equation*}
Recently, Quahhabi and Tahri studied generalized Heisenberg algebra and quantum operator functions. The generalized Heisenberg algebra consists of the algebra $GHA_{f}$ generated by $\{a,\ap,H\}$ satisfying the following commutation relations: 
\begin{equation}
aH=f(H)a,\quad H\ap=\ap f(H),	\label{6}
\end{equation}
and 
\begin{equation}
[a,\ap]=f(H)-H, \label{7}	
\end{equation}
where $H$ is the self-conjugate Hamiltonian of the system and $f(H)$ is an analytic function of $H$, called the characteristic function of the algebra (see [11]). Let $f(H)=H+1$. Then we have 
\begin{equation}
aH=Ha+a,\quad H\ap=\ap H+\ap,\quad [a,\ap]=1,\label{8}
\end{equation}
so that $\text{GHA}_f$ reduces to the classical Heisenberg-Weyl algebra.
For any positive integer $k$, we define the power composition of the function $f$ by 
\begin{equation}
f^{k}=\underbrace{f\circ f\circ \cdots \circ f}_{k-\mathrm{times}},\ \mathrm{and}\ f^{-k}= \underbrace{f^{-1}\circ f^{-1}\circ \cdots \circ f^{-1}}_{k-\mathrm{times}},\quad (\mathrm{see}\ [11]), \label{9}
\end{equation}
where $f^{-1}$ is the compositional inverse of $f$ and $f^{0}$ denotes the identity function. \par 
For $k\in\mathbb{Z}$, quantum operator functions are defined by 
\begin{equation}
[k]_{f}(H)=f^{k}(H)-H,\quad (\mathrm{see}\ [5,8,11]).\label{10}
\end{equation}
Note that 
\begin{equation*}
[-k]_{f}(H)=[k]_{f^{-1}}(H). 
\end{equation*}
When $f(H)=H+1$, we have $[k]_{f}(H)=k$. \par 
In addition, quantum operator functions of the second kind are given by 
\begin{equation}
(k)_{f}(H)=H-f^{-k}(H)=[k]_{f}\big(f^{-k}(H)\big),\quad (\mathrm{see}\ [8,11]). \label{11}
\end{equation}
It is known that 
\begin{align}
&a^{k}[n]_{f}(H)=[n]_{f}\big(f^{k}(H)\big)a^{k},\label{12}\\
&[n]_{f}(H)\big(\ap\big)^{k}=\big(\ap\big)^{k}[n]_{f}\big(f^{k}(H)\big), \nonumber\\
&\big[a,(\ap)^{k}\big]=\big(\ap\big)^{k-1}[k]_{f}(H),\quad \big[a^{k},\ap\big]=[k]_{f}(H)a^{k-1},\quad (\mathrm{see}\ [11]). \nonumber	
\end{align}
Fock space representation of the generalized Heiseberg algebra is obtained with respect to the orthonormal eigenstate basis $\{|n\rangle\}_{n\in\mathbb{N}\cup\{0\}}$, of the Hamiltonian $H$. That is 
\begin{equation}
H|m\rangle=\alpha_{m}|m\rangle,\quad (\mathrm{see} [11]).\label{13}	
\end{equation}
The whole space is constructed from the lowest eigenvalue of $H$ corresponding to the vacuum state $|0\rangle$. \par 
On the one hand, the characteristic function generates all the eigenvalues of $H$ by the successive iteration from $\alpha_{0}$. Indeed, from \eqref{6} and \eqref{13}, we note that 
\begin{equation}
\alpha_{m}=f(\alpha_{m-1})=f^{m}(\alpha_{0}).\label{14}
\end{equation}
On the other hand, the operators $\ap$ and $a$ generate all the eigenvalues of $H$ by successive actions as in the following:
\begin{align}
&a|m \rangle=N_{m-1}|m-1 \rangle, \quad a|0 \rangle=0 \quad(N_{-1}=0), \label{15} \\
& \ap|m \rangle=N_{m}|m+1 \rangle. \nonumber
\end{align}
We observe from \eqref{15} that
\begin{equation}
\fa|m \rangle=\ap N_{m-1}|m-1 \rangle=N_{m-1}\ap|m-1 \rangle=N_{m-1}^{2}|m \rangle. \label{16}
\end{equation}
Similarly, we show from \eqref{15} that
\begin{equation}
a\ap|m \rangle=N_{m}^{2}|m \rangle . \label{17}
\end{equation}
Noting $\fa=a\ap-f(H)+H$ from \eqref{7}, and using \eqref{13} and \eqref{17}, we have
\begin{equation}
\fa|m \rangle=a\ap-f(H)+H|m \rangle=N_{m}^{2}-f(\alpha_{m})+\alpha_{m}|m \rangle. \label{18}
\end{equation}
Thus, from \eqref{14}, \eqref{16} and \eqref{18}, we obtain
\begin{equation}
N_{m}^{2}-f^{m+1}(\alpha_{0})=N_{m-1}^{2}-f^{m}(\alpha_{0}). \label{19}
\end{equation}
We now deduce from \eqref{19} that
\begin{equation}
N_{m-1}^{2}-f^{m}(\alpha_{0})=N_{-1}^{2}-f^{0}(\alpha_{0})=-\alpha_{0}. \label{20}
\end{equation}
Therefore, we have shown that 
\begin{align}
N_{m-1}^{2}&=f^{m}(\alpha_{0})-\alpha_{0}=\alpha_{m}-\alpha_{0} \label{21} \\
&=[m]_{f}(\alpha_{0})=\big\langle 0\big|[m]_{f}(H)\big|0\big\rangle,\quad (m\ge 0). \nonumber
\end{align}
From \eqref{15} and \eqref{21}, the actions of $a$ and $\ap$ on number states in Fock space are given by 
\begin{align}
&a|m\rangle=\sqrt{[m]_{f}(\alpha_{0})}\big|m-1\big\rangle,\quad a|0\rangle=0,	\label{22} \\
&\ap|m\rangle=\sqrt{[m+1]_{f}(\alpha_{0})}\big|m+1\big\rangle,\quad (\mathrm{see}\ [11]).	\nonumber
\end{align}
Let $N_{f}=\fa$. Then, by \eqref{16} and \eqref{21}, we have
\begin{equation}
N_{f}|m\rangle=[m]_{f}(\alpha_{0})\big|m\big\rangle,\quad (\mathrm{see}\ [11]). \label{23}
\end{equation}
From \eqref{13}, \eqref{21} and \eqref{23}, we have
\begin{equation}
H|m \rangle=\alpha_{m}|m \rangle=\alpha_{0}+[m]_{f}(\alpha_{0})|m \rangle 
=\alpha_{0}+N_{f}|m \rangle. \label{24}
\end{equation}
By \eqref{24}, we get 
\begin{equation}
H=N_{f}+\alpha_{0}. \label{25}
\end{equation}
Recently, Ouahhabi and Tahri introduced the generalized Stirling operators of the second kind defined by \begin{equation}
 N_{f}^{n}=\big(\fa\big)^{n}=\sum_{k=0}^{n}{n \brace k}_{f}(H)\big(\ap\big)^{k}a^{k},\quad (n\ge 0),\quad (\mathrm{see}\ [11]). \label{26}	
 \end{equation}
When $f(H)=H+1$, we have 
\begin{equation*}
N^{n}=\big(\fa\big)^{n}=\sum_{k=0}^{n}{n \brace k}\big(\ap\big)^{k}a^{k},\quad (\mathrm{see}\ [11]). 
\end{equation*}
From \eqref{26}, we note that 
\begin{equation}
{n \brace 0}_{f}(H)=\delta_{n,0},\quad {n \brace k}_{f}(H)=0,\quad\textrm{if $k>n$}, \label{27}
\end{equation}
where $\delta_{n,k}$ is the Kronecker delta (see [11]). \par 

\section{Generalized Bell polynomial operators arising from generalized normal ordering}
We define the {\it{quantum operator factorials}} by 
\begin{equation}
[n]_{f}(\alpha_{0})!=[n]_{f}(\alpha_{0})[n-1]_{f}(\alpha_{0})\cdots[1]_{f}(\alpha_{0}),\quad (n \ge 1),\quad [0]_{f}(\alpha_{0})!=1. \label{28}
\end{equation}
Then, we introduce the {\it{quantum operator exponential function}} defined by 
\begin{equation}
e_{f}^{x}=\sum_{n=0}^{\infty}\frac{x^{n}}{[n]_{f}(\alpha_{0})!}. \label{29}	
\end{equation}
We recall the coherent states:
\begin{equation*}
a|z\rangle=z|z \rangle,\quad \textrm{equivalently}\quad \langle z|\ap=\langle z|\overline{z},\quad\textrm{and $\langle z|z\rangle=1$, $(z\in\mathbb{C})$.}
\end{equation*}
For the coherent state $|z\rangle$, we assume that 
\begin{equation}
|z\rangle=\sum_{n=0}^{\infty}q_{n}|n\rangle.\label{30}
\end{equation}
Thus, by \eqref{30}, we get 
\begin{equation}
a|z\rangle=\sum_{n=0}^{\infty}q_{n}a|n\rangle=\sum_{n=0}^{\infty}\sqrt{[n]_{f}(\alpha_{0})}q_{n}\big|n-1\big\rangle,\label{31}	
\end{equation}
and 
\begin{equation}
a|z\rangle=z|z\rangle=z\sum_{n=0}^{\infty}q_{n}|n\rangle=\sum_{n=0}^{\infty}zq_{n-1}|n-1\rangle.\label{32}
 \end{equation}
From \eqref{31} and \eqref{32}, we have 
\begin{align}
q_{n}&=\frac{z}{\sqrt{[n]_{f}(\alpha_{0})}}q_{n-1}=\frac{z^{2}}{\sqrt{[n]_{f}(\alpha_{0})}\sqrt{[n-1]_{f}(\alpha_{0})}}q_{n-2}=\cdots \label{33}\\
&=\frac{z^{n}}{ \sqrt{[n]_{f}(\alpha_{0})}\sqrt{[n-1]_{f}(\alpha_{0})}\cdots\sqrt{[1]_{f}(\alpha_{0})}}q_{0}=q_{0}\frac{z^{n}}{\sqrt{[n]_{f}(\alpha_{0})!}}.\nonumber
\end{align}
By \eqref{30} and \eqref{33}, we get 
\begin{equation}
|z\rangle=q_{0}\sum_{n=0}^{\infty}\frac{z^{n}}{\sqrt{[n]_{f}(\alpha_{0})!}}\bigg|n\bigg\rangle.\label{34}
\end{equation}
From \eqref{34}, we note that 
\begin{align}
1=\langle z|z\rangle &=\overline{q_{0}}\sum_{m=0}^{\infty}\frac{(\overline{z})^{m}}{\sqrt{[m]_{f}(\alpha_{0})!}}q_{0}\sum_{n=0}^{\infty}\frac{z^{n}}{\sqrt{[n]_{f}(\alpha_{0})!}}\langle m|n\rangle \label{35}\\
&=|q_{0}|^{2}\sum_{n=0}^{\infty}\frac{|z|^{2n}}{[n]_{f}(\alpha_{0})!}=|q_{0}|^{2}e_{f}^{|z|^{2}}. \nonumber
\end{align}
Thus, by \eqref{35}, we get 
\begin{equation}
|q_{0}|=e_{f}^{-\frac{|z|^{2}}{2}}.\label{36}
\end{equation}
From \eqref{34} and \eqref{36}, we have 
\begin{equation}
|z\rangle=e_{f}^{-\frac{|z|^{2}}{2}}\sum_{n=0}^{\infty}\frac{z^{n}}{\sqrt{[n]_{f}(\alpha_{0})!}}\bigg|n\bigg\rangle. \label{37}	
\end{equation}
Therefore, by \eqref{37}, we obtain the following theorem. 
\begin{theorem}
For $z\in\mathbb{C}$, we have 
\begin{equation*}
|z\rangle=e_{f}^{-\frac{|z|^{2}}{2}}\sum_{n=0}^{\infty}\frac{z^{n}}{\sqrt{[n]_{f}(\alpha_{0})!}}\bigg|n\bigg\rangle.
\end{equation*}
\end{theorem}
From \eqref{26}, we have 
\begin{equation*}
N_{f}^{n}=\big(\fa\big)^{n}=\sum_{k=0}^{n}{n \brace k}_{f}(H)(\ap)^{k}a^{k}. 
\end{equation*}
Thus, we see that 
\begin{align}
N_{f}^{n}|m\big\rangle &=\sum_{l=0}^{n}{n \brace l}_{f}(H)\big(\ap\big)^{l}a^{l}\big|m\big\rangle \label{38}\\
&=\sum_{l=0}^{n}{n \brace l}_{f}(\alpha_{m})\big([m]_{f}(\alpha_{0})\big)_{l}\big|m\big\rangle,\nonumber
\end{align}
where the generalized falling factorial operators are given by 
\begin{equation}
\begin{aligned}
&\big([m]_{f}(\alpha_{0})\big)_{0}=1,\quad \\
&\big([m]_{f}(\alpha_{0})\big)_{l}=[m]_{f}(\alpha_{0})[m-1]_{f}(\alpha_{0})\cdots [m-l+1]_{f}(\alpha_{0}),\ (l\ge 1). 
\end{aligned}
\label{39}
\end{equation}
On the other hand, by \eqref{23}, we get 
\begin{align}
N_{f}^{n}|m\rangle &=\big(\fa\big)^{n}\big|m\big\rangle=\underbrace{(\fa)(\fa)\cdots(\fa)}_{n-\mathrm{times}}\big|m\big\rangle \label{40}\\
&=\underbrace{[m]_{f}(\alpha_{0})\cdot [m]_{f}(\alpha_{0})\cdots [m]_{f}}_{n-\mathrm{times}}(\alpha_{0})\big|m\big\rangle=\big([m]_{f}(\alpha_{0})\big)^{n}\big|m\big\rangle. \nonumber
\end{align}
Therefore, by \eqref{38} and \eqref{40}, we obtain the following theorem. 
\begin{theorem}
For $n\ge 0$, we have 
\begin{equation*}
\Big([m]_{f}(\alpha_{0})\Big)^{n}=\sum_{l=0}^{n}{n \brace l}_{f}(\alpha_{m})\big([m]_{f}(\alpha_{0})\big)_{l}. 
\end{equation*}
\end{theorem}
From the definition of coherent states $|z\rangle,\ (z \in \mathbb{C})$, we note that 
\begin{align}
\big\langle z\big|\big(\fa\big)^{n}\big|z\big\rangle &=\sum_{l=0}^{n}{n \brace l}_{f}(\alpha_{0}+|z|^{2})\big\langle z\big|(\ap)^{l}a^{l}\big|z\big\rangle \label{41}\\
&=\sum_{l=0}^{n}{n \brace l}_{f}(\alpha_{0}+|z|^{2})\big(\overline{z}\big)^{l}z^{l}\langle z|z\rangle \nonumber\\
&=\sum_{l=0}^{n}{n \brace l}_{f}(\alpha_{0}+|z|^{2})|z|^{2l},\quad (n\ge 0).\nonumber
\end{align}
Therefore, by \eqref{41}, we obtain the following theorem. 
\begin{theorem}
For $n\ge 0$, we have 
\begin{equation*}
\big\langle z\big|\big(\fa\big)^{n}\big|z\big\rangle=\sum_{l=0}^{n}{n \brace l}_{f}(\alpha_{0}+|z|^{2})|z|^{2l}. 
\end{equation*}
\end{theorem}
In view of \eqref{2}, we define the {\it{generalized Bell polynomial operators arising from generalized normal ordering}} as
\begin{equation}
\phi_{n,f}(x)(H)=\sum_{l=0}^{n}{n \brace l}_{f}(H)x^{l},\quad (n\ge 0). \label{42}	
\end{equation}
When $x=1,\ \phi_{n,f}(H)=\phi_{n,f}(1)(H)$ are called the {\it{generalized Bell number operators arising from generalized normal ordering}}. \par 
From \eqref{41} and \eqref{42}, we obtain the following theorem. 
\begin{theorem}
For $n\ge 0$, we have 
\begin{equation}
\big\langle z\big|\big(\fa\big)^{n}\big|z\big\rangle=\langle z|N_{f}^{n}|z\rangle=\phi_{n,f}(|z|^{2})(\alpha_{0}+|z|^{2}).\label{43}
\end{equation}
\end{theorem}
Let $f(t,z)=\big\langle z\big|e^{\fa t}\big|z\big\rangle$. Then we have 
\begin{equation}
f(t,z)=\sum_{k=0}^{\infty}\frac{t^{k}}{k!}\big\langle z\big|\big(\fa\big)^{k}\big|z\big\rangle=\sum_{k=0}^{\infty}\frac{t^{k}}{k!}\phi_{k,f}(|z|^{2})(\alpha_{0}+|z|^{2}). \label{44}
\end{equation}
From \eqref{43}, we note that 
\begin{align}
\phi_{k,f}(|z|^{2})(H)&=\big\langle z\big|N_{f}^{k}\big|z\big\rangle \label{45}	\\
&=e_{f}^{-\frac{|z|^{2}}{2}}\cdot e_{f}^{-\frac{|z|^{2}}{2}}\sum_{m,n=0}^{\infty}\frac{(\overline{z})^{m}z^{n}}{\sqrt{[m]_{f}(\alpha_{0})!} \sqrt{[n]_{f}(\alpha_{0})!} }\big\langle m\big|\big(\fa\big)^{k}\big|n\big\rangle \nonumber\\
&=e_{f}^{-|z|^{2}}\sum_{m,n=0}^{\infty} \frac{(\overline{z})^{m}z^{n}}{\sqrt{[m]_{f}(\alpha_{0})!} \sqrt{[n]_{f}(\alpha_{0})!} }\big([n]_{f}(\alpha_{0})\big)^{k}\langle m|n\rangle \nonumber\\
&=e_{f}^{-|z|^{2}}\sum_{n=0}^{\infty}\frac{|z|^{2n}}{[n]_{f}(\alpha_{0})!}\big([n]_{f}(\alpha_{0})\big)^{k}. \nonumber
\end{align}
Therefore, by \eqref{45}, we obtain the following theorem. 
\begin{theorem}
For $k\ge 0,\ z\in\mathbb{C}$, we have 
\begin{equation}
\phi_{k,f}\big(|z|^{2}\big)(\alpha_{0}+|z|^{2})=e_{f}^{-|z|^{2}}\sum_{n=0}^{\infty}\frac{|z|^{2n}}{[n]_{f}(\alpha_{0})!}\big([n]_{f}(\alpha_{0})\big)^{k}. \label{46}
\end{equation}
In particular, for $|z|=1$, 
\begin{equation*}
\phi_{k,f}(\alpha_{0}+1)=e_{f}^{-1}\sum_{n=0}^{\infty}\frac{1}{[n]_{f}(\alpha_{0})!}\big([n]_{f}(\alpha_{0})\big)^{k}. 
\end{equation*}
\end{theorem}
When $f(H)=H+1$, we note that 
\begin{equation}
[n]_{f}(\alpha_{0})=n,\quad e_{f}^{-|z|^{2}}=e^{-|z|^{2}}.\label{47}	
\end{equation}
By \eqref{46} and \eqref{47}, we get (see \eqref{2})
\begin{equation}
\phi_{k}\big(|z|^{2}\big)=e^{-|z|^{2}}\sum_{n=0}^{\infty}\frac{|z|^{2n}}{n!}n^{k},\ (k\ge 0). \label{48}	
\end{equation}

\section{Conclusion}
In this paper, we have successfully extended the classic combinatorial structures of quantum optics into the framework of the generalized Heisenberg algebra $\text{GHA}_f$. By evaluating the normal ordering of the deformed number operator $N_f^n$, we established a rigorous algebraic foundation for generalized Stirling operators of the second kind and generalized Bell polynomial operators. Our main results culminated in the exact derivation of generalized coherent states and their structural identities. We demonstrated that the expectation value of the $n$-th power of the generalized number operator within these coherent states is governed precisely by the generalized Bell polynomial operators $\phi_{n,f}(|z|^2)(H)$. These findings clarify the algebraic behavior of systems with non-linear spectra and provide a robust operator calculus framework. Future work will focus on applying these generalized combinatorial operators to evaluate the non-classical properties, such as squeezing and antibunching, of specific physical systems modeled by $\text{GHA}_f$.

\end{document}